\documentclass[a4paper,reqno,12pt,oneside]{amsart}
\usepackage{amsmath,amsrefs,amssymb}

\numberwithin{equation}{section}

\DeclareMathOperator{\divergence}{div}
\DeclareMathOperator{\dist}{dist}
\newcommand{\OO}{\text{O}}
\newcommand{\oo}{\text{o}}
\newcommand{\N}{\mathbb{N}}
\newcommand{\R}{\mathbb{R}}
\renewcommand{\S}{\mathbb{S}}
\newcommand{\<}{\left<}
\renewcommand{\>}{\right>}
\renewcommand{\[}{\left[}
\renewcommand{\]}{\right]}
\renewcommand{\(}{\left(}
\renewcommand{\)}{\right)}

\newtheorem{theorem}{Theorem}[section]
\newtheorem{step}[theorem]{Step}
\newtheorem{lemma}[theorem]{Lemma}

\begin{document}

\title[Results on the Guan--Li mean curvature flow]{Convergence result and blow-up examples for the Guan--Li mean curvature flow on warped product spaces}

\author{J\'er\^ome V\'etois}

\address{J\'er\^ome V\'etois, McGill University Department of Mathematics and Statistics, 805 Sherbrooke Street West, Montreal, Quebec H3A 0B9, Canada.}
\email{jerome.vetois@mcgill.ca}

\date{July 28, 2018.}

\thanks{To appear in {\it Communications in Analysis and Geometry}.}

\begin{abstract}
We examine the question of convergence of solutions to a geometric flow which was introduced by Guan and Li~\cite{GuanLi} for starshaped hypersurfaces in space forms and generalized by Guan, Li, and Wang~\cite{GuanLiWang} to the case of warped product spaces. We obtain a convergence result under a condition on the optimal modulus of continuity of the initial data. Moreover we show by examples that this condition is optimal at least in the one-dimensional case.
\end{abstract}

\maketitle

\section{Introduction and main results}\label{Sec1}\addtolength{\textheight}{-6pt}

Let $n\ge1$, $\(\S^n,g_{\S^n}\)$ be the standard $n$--sphere, $I\subset\R$ be a closed interval, and $\(N,\overline{g}\)$ be the warped product of $\S^n$ and $I$ equipped with
$$\overline{g}=\phi\(\rho\)^2g_{\S^n}+d\rho^2,$$
where $\phi:I\to\(0,\infty\)$ is a smooth warping function. We consider the following flow which was introduced by Guan and Li~\cite{GuanLi} in the case of space forms and generalized by Guan, Li, and Wang~\cite{GuanLiWang} to the case of warped product spaces (see also Cant~\cite{Cant} in case $n=1$):
$$\partial_tF=\(n\phi'-Hu\)\nu,$$
where $T\in\(0,\infty\]$, $\(F\(\cdot,t\)\)_{t\in\[0,T\)}$ is a smooth family of embeddings into $N$ which defines smooth hypersurfaces $\(M_t\)_{t\in\[0,T\)}$ and $H$, $u$, $\nu$ are the mean curvature, support function, and outward unit normal vector field, respectively, of the hypersurfaces $\(M_t\)_{t\in\[0,T\)}$. A crucial property of this flow is that it preserves the volume enclosed by the initial hypersurface while monotonically decreasing the area (see~\cite{GuanLi}*{Proposition~3.5}). 

\smallskip
Throughout this paper, we assume that the hypersurfaces $\(M_t\)_{t\in\[0,T\)}$ are starshaped, i.e. for every $t\in\[0,T\)$, $M_t$ is the graph of a function $\rho\(\cdot,t\):\S^n\to\R$. We then obtain (see the formulas in~\cite{GuanLiWang}*{Section~3}) that $\rho$ solves the initial value problem
\begin{equation}\label{Eq1}
\left\{\begin{aligned}&\partial_t\rho=\divergence\(\frac{\nabla\rho}{\sqrt{\phi\(\rho\)^2+\left|\nabla\rho\right|^2}}\)+n\frac{\phi'\(\rho\)}{\phi\(\rho\)}\frac{\left|\nabla\rho\right|^2}{\sqrt{\phi\(\rho\)^2+\left|\nabla\rho\right|^2}}&&\text{in }D_T\\
&\rho\(\cdot,0\)=\rho_0&&\text{on }\S^n,
\end{aligned}\right.
\end{equation}
where $\divergence:=\divergence_{g_{\S^n}}$, $\left|\cdot\right|=\left|\cdot\right|_{g_{\S^n}}$, $D_T:=\S^n\times\(0,T\)\to\R$, and $\rho_0$ is the radial function of $M_0$. It follows from classical theory of parabolic equations that for every $\rho_0\in C^\infty\(\S^n,I\)$, there exists a unique solution $\rho\in C^\infty\(\overline{D_T}\)$ of \eqref{Eq1} for small $T>0$. Moreover, a straightforward application of the maximum principle gives $\rho\(D_T\)\subseteq I$.

\smallskip
The following result has been obtained by Guan and Li~\cite{GuanLi}  in the case of space forms and generalized by Guan, Li, and Wang~\cite{GuanLiWang} to the case of warped product spaces:

\begin{theorem}\label{Th0} (Guan and Li~\cite{GuanLi}, Guan, Li, and Wang~\cite{GuanLiWang})
Let $I\subset\R$ be a closed interval, $\phi\in C^\infty\(I,\(0,\infty\)\)$, and $n\ge1$. Assume that
\begin{equation}\label{Eq2}
\phi'^2-\phi\phi''\ge0\quad\text{in }I.
\end{equation}
Then for any $\rho_0\in C^\infty\(\S^n,I\)$, the solution of \eqref{Eq1} exists for all time and converges exponentially to a constant i.e. $T=\infty$ and there exist $\rho_\infty\in I$, $C,\eta>0$ such that $\left|\rho\(x,t\)-\rho_\infty\right|\le Ce^{-\eta t}$ for all $\(x,t\)\in D_\infty$.
\end{theorem}

This result has been successfully used in~\cites{Cant,GuanLi,GuanLiWang} to solve isoperimetric problems in warped product spaces. As is explained in~\cite{GuanLiWang}*{Proposition~6.1}, the condition \eqref{Eq2} is strongly related to the notion of photon sphere in general relativity.

\smallskip
In this paper, we investigate the case where the condition \eqref{Eq2} is not satisfied. In this case, we obtain a convergence result under a barrier condition on the optimal modulus of continuity of $\rho_0$, namely
$$\omega_{\rho_0}\(\theta\):=\sup\left\{\left|\rho_0\(y\)-\rho_0\(x\)\right|\,:\,x,y\in\S^n\text{ and }\dist_{\S^n}\(x,y\)=\theta\right\}$$
for all $\theta\in\[0,\pi\]$. Here $\dist_{\S^n}$ denotes the distance on $\S^n$ with respect to the standard metric. We obtain the following result:

\begin{theorem}\label{Th1}
Let $I\subset\R$ be a closed interval, $\phi\in C^\infty\(I,\(0,\infty\)\)$, and $n\ge1$. Then there exists $\lambda_0>0$ such that for any $\rho_0\in C^\infty\(\S^n,I\)$, if
\begin{equation}\label{Eq3}
\omega_{\rho_0}\(\theta\)\le\lambda_0\theta^{1/2}\qquad\forall\theta\in\[0,\pi\],
\end{equation} 
then the solution of \eqref{Eq1} exists for all time and converges exponentially to a constant.
\end{theorem}

We prove Theorem~\ref{Th1} in Section~\ref{Sec2} by using an approach based on Kruzhkov's doubling variable technique \cite{Kru} and inspired by the works of Andrews and Clutterbuck~\cites{AndClu1,AndClu2,AndClu3,AndClu4}. As in the papers of Cant~\cites{Cant}, Guan and Li~\cite{GuanLi}, and Guan, Li, and Wang~\cite{GuanLiWang}, Theorem~\ref{Th1} can be applied to solve isoperimetric problems in the warped product space $\(N,\overline{g}\)$ provided $\phi'^2-\phi\phi''\le1$ in $I$, which is a necessary condition for the isoperimetric inequality (see Li and Wang~\cites{LiWang}).

\smallskip
The following result, obtained in case $n=1$, shows the optimality of the exponent $1/2$ in \eqref{Eq3}:

\begin{theorem}\label{Th2}
Assume that $n=1$, $0\in I$, $\phi$ is even, and $\phi''\(0\)>0$. Then for any $\sigma\in\(0,1/2\)$  $\lambda>0$, there exist $\rho_0\in C^\infty\(\S^n,I\)$ such that 
\begin{equation}\label{Eq4}
\omega_{\rho_0}\(\theta\)\le\lambda\theta^\sigma\qquad\forall\theta\in\[0,\pi\]
\end{equation} 
and the solution of \eqref{Eq1} is such that $\partial_x\rho$ blows up in finite time i.e. $\sup_{D_t}\left|\partial_x\rho\right|\to\infty$ as $t\to T$ for some $T\in\(0,\infty\)$.
\end{theorem}

We prove Theorem~\ref{Th2} in Section~\ref{Sec3}. As far as the author knows, this is the first existence result of blowing-up solutions for \eqref{Eq1}. The high nonlinearity of the flow makes it difficult to construct examples of blowing-up solutions. Here, the solutions that we construct are periodic, with a large number of oscillations. Our existence result relies on the construction of a suitable family of barrier functions on a small arc of $\S^1$ with zero boundary condition. We then exploit the symmetry of the warping function to extend our solutions to the whole $\S^1$.

\smallskip\noindent
{\bf Acknowledgments.} The author is very grateful to Pengfei Guan for many enlightening discussions and helpful advice during the preparation of this paper.

\section{Proof of Theorem~\ref{Th1}}\label{Sec2}

This section is devoted to the proof of Theorem~\ref{Th1}. As in the paper of Guan and Li~\cite{GuanLi}, it will be convenient to use the change of functions
\begin{equation}\label{Th1Eq2}
\gamma=\Gamma\(\rho\):=\int_{\overline\rho}^\rho\frac{ds}{\phi\(s\)}\qquad\text{and}\qquad\psi:=\phi\circ\Gamma^{-1},
\end{equation}
where $\overline\rho\in I$ is fixed. By differentiating, we obtain $\nabla\rho=\psi\(\gamma\)\nabla\gamma$ and $\phi'\(\rho\)=\psi'\(\gamma\)/\psi\(\gamma\)$. Hence the problem \eqref{Eq1} becomes
\begin{equation}\label{Th1Eq3}
\left\{\begin{aligned}&\partial_t\gamma=\frac{1}{\psi\(\gamma\)}\divergence\(\frac{\nabla\gamma}{\sqrt{1+\left|\nabla\gamma\right|^2}}\)+n\frac{\psi'\(\gamma\)}{\psi\(\gamma\)^2}\frac{\left|\nabla\gamma\right|^2}{\sqrt{1+\left|\nabla\gamma\right|^2}}&&\text{in }D_T\\
&\gamma\(\cdot,0\)=\gamma_0&&\text{on }\S^n,
\end{aligned}\right.
\end{equation}
where 
$$\gamma_0:=\int_{\overline\rho}^{\rho_0}\frac{ds}{\phi\(s\)}\,.$$
A straightforward application of the maximum principle gives $\gamma\(D_T\)\subseteq\gamma_0\(\S^n\)$ i.e.
$$\min_{\S^n}\gamma_0\le\gamma\(x,t\)\le\max_{\S^n}\gamma_0\qquad\forall\(x,t\)\in D_T.$$
We assume that $\rho_0\(\S^n\)\subseteq I$ and we let $\lambda>0$ be such that 
\begin{equation}\label{Th1Eq4}
\omega_{\rho_0}\(\theta\)\le\lambda\sqrt\theta\qquad\forall\theta\in\[0,\pi\].
\end{equation}
Since $\rho_0\in C^\infty\(\S^n,I\)$, we obtain that there exists $\Lambda>0$ such that
\begin{equation}\label{Th1Eq4a}
\omega_{\rho_0}\(\theta\)\le\Lambda\theta\qquad\forall\theta\in\[0,\pi\].
\end{equation}
For every $\delta\in\(0,\lambda^2/\Lambda^2\)$, an easy study of functions gives that
\begin{equation}\label{Th1Eq4b}
\Lambda\theta\le2\lambda\big(\sqrt{\delta+\theta}-\sqrt\delta\big)\qquad\forall\theta\in\Big(0,\frac{4\lambda}{\Lambda^2}\big(\lambda-\Lambda\sqrt\delta\big)\Big)
\end{equation} 
and
\begin{equation}\label{Th1Eq4c}
\sqrt\theta\le2\big(\sqrt{\delta+\theta}-\sqrt\delta\big)\qquad\forall\theta>\frac{16}{9}\delta.
\end{equation} 
It follows from \eqref{Th1Eq4}--\eqref{Th1Eq4c} that we can choose $\delta\in\(0,1\)$ small enough such that
\begin{equation}\label{Th1Eq5}
\omega_{\rho_0}\(\theta\)\le2\lambda\big(\sqrt{\delta+\theta}-\sqrt\delta\big)\qquad\forall\theta\in\[0,\pi\].
\end{equation}
By using the mean value theorem, it follows from \eqref{Th1Eq5} that
\begin{equation}\label{Th1Eq6}
\omega_{\gamma_0}\(\theta\)\le2\overline\lambda\big(\sqrt{\delta+\theta}-\sqrt\delta\big)\quad\forall \theta\in\[0,\pi\],
\end{equation}
where 
$$\overline\lambda:=\lambda\sup_I\frac{1}{\phi}.$$
We will show that if $\lambda$ is smaller than a constant $\lambda_0$ depending only on $I$, $\phi$, and $n$, then $\left|\nabla\gamma\right|$ is bounded above by an exponentially decaying function. We will use an approach based on Kruzhkov's doubling variable technique \cite{Kru}. This approach was successfully used in the works of Andrews and Clutterbuck~\cites{AndClu1,AndClu2,AndClu3,AndClu4} to obtain sharp estimates on the gradient and modulus of continuity of solutions to quasilinear parabolic equations. We fix $\eta>0$ and we define
$$\kappa\(\theta,t\):=2\overline\lambda\big(\sqrt{\delta+\theta}-\sqrt\delta\big)e^{-\eta t}\qquad\forall\(\theta,t\)\in\[0,\pi\]\times\[0,T\],$$
where $\delta$ and $\overline\lambda$ are as above. We then define
$$Z\(x,y,t\):=\gamma\(y,t\)-\gamma\(x,t\)-\kappa\(d\(x,y\),t\)\qquad\forall\(x,y,t\)\in U_T,$$
where $U_T:=\(\S^n\)^2\times\[0,T\]$. It follows from \eqref{Th1Eq6} that $Z\(x,y,0\)\le0$ for all $x,y\in\S^n$. In what follows, we will show that if $\eta$ and $\lambda$ are small enough, then $Z\(x,y,t\)\le0$ for all $\(x,y,t\)\in U_T$. We assume by contradiction that $Z$ is not everywhere nonpositive in $U_T$. Then we obtain that for small $\varepsilon>0$, there exists $\(x_\varepsilon,y_\varepsilon,t_\varepsilon\)\in U_T$ such that
\begin{equation}\label{Th1Eq7}
Z\(x_\varepsilon,y_\varepsilon,t_\varepsilon\)=\max\(\left\{Z\(x,y,t\):\,\,x,y\in\S^n\,\,\text{and}\,\,t\le t_\varepsilon\right\}\)=\varepsilon.
\end{equation}
We define $\theta_\varepsilon:=d\(x_\varepsilon,y_\varepsilon\)$. 

\smallskip
As a first step, we obtain the following result:

\begin{step}\label{Th1St1}
$\varepsilon=\oo\(\theta_\varepsilon\)$ as $\varepsilon\to0$.
\end{step}

\proof[Proof of Step~\ref{Th1St1}]
Assume by contradiction that there exists a sequence $\(\varepsilon_\alpha\)_{\alpha\in\mathbb{N}}$ such that $\varepsilon_\alpha>0$, $\theta_{\varepsilon_\alpha}=\OO\(\varepsilon_\alpha\)$, and $\varepsilon_\alpha\to0$ as $\alpha\to\infty$. Since $\kappa\(0,t_{\varepsilon_\alpha}\)=0$, by applying the mean value theorem, we obtain that there exist $\zeta_\alpha,\xi_\alpha\in\(0,\theta_{\varepsilon_\alpha}\)$ such that 
\begin{equation}\label{Th1St1Eq1}
\kappa\(\theta_{\varepsilon_\alpha},t_{\varepsilon_\alpha}\)=\partial_\theta\kappa\(\zeta_\alpha,t_{\varepsilon_\alpha}\)\theta_{\varepsilon_\alpha}
\end{equation}
and
\begin{equation}\label{Th1St1Eq2}
\gamma\(y_{\varepsilon_\alpha},t_{\varepsilon_\alpha}\)-\gamma\(x_{\varepsilon_\alpha},t_{\varepsilon_\alpha}\)=\<\nabla_x\gamma\(\tau_{\varepsilon_\alpha}\(\xi_\alpha\),t_{\varepsilon_\alpha}\),\tau'_{\varepsilon_\alpha}\(\xi_\alpha\)\>\theta_{\varepsilon_\alpha},
\end{equation}
where $\tau_{\varepsilon_\alpha}:\[0,\theta_{\varepsilon_\alpha}\]\to\S^n$ is a minimizing geodesic from $x_{\varepsilon_\alpha}$ to $y_{\varepsilon_\alpha}$. It follows from \eqref{Th1Eq7}, \eqref{Th1St1Eq1}, \eqref{Th1St1Eq2}, and Cauchy--Schwartz inequality that $\theta_{\varepsilon_\alpha}\ne0$ and
\begin{equation}\label{Th1St1Eq3}
\frac{\varepsilon_\alpha}{\theta_{\varepsilon_\alpha}}\le\left|\nabla_x\gamma\(\tau_{\varepsilon_\alpha}\(\xi_\alpha\),t_{\varepsilon_\alpha}\)\right|-\partial_\theta\kappa\(\zeta_\alpha,t_{\varepsilon_\alpha}\).
\end{equation}
Since $\(t_{\varepsilon_\alpha}\)_{\alpha\in\mathbb{N}}$ is decreasing, we obtain $t_{\varepsilon_\alpha}\to t_0$ for some $t_0\ge0$. Since $\theta_{\varepsilon_\alpha}\to0$, we obtain $\zeta_\alpha,\xi_\alpha\to0$. Moreover up to a subsequence $x_{\varepsilon_\alpha},y_{\varepsilon_\alpha}\to x_0\in\S^n$. By passing to the limit into \eqref{Th1St1Eq3}, we then obtain 
\begin{equation}\label{Th1St1Eq4}
\limsup_{\alpha\to\infty}\frac{\varepsilon_\alpha}{\theta_{\varepsilon_\alpha}}\le\left|\nabla_x\gamma\(x_0,t_0\)\right|-\partial_\theta\kappa\(0,t_0\).
\end{equation}
On the other hand, by passing to the limit into \eqref{Th1Eq7}, first as $\varepsilon\to0$ and  then as $x,y\to x_0$, we obtain
\begin{equation}\label{Th1St1Eq5}
\left|\nabla_x\gamma\(x_0,t_0\)\right|\le\partial_\theta\kappa\(0,t_0\).
\end{equation}
By putting together \eqref{Th1St1Eq4} and \eqref{Th1St1Eq5}, we obtain a contradiction with $\theta_{\varepsilon_\alpha}=\OO\(\varepsilon_\alpha\)$. This ends the proof of Step~\ref{Th1St1}.
\endproof

We then prove the following result:

\begin{step}\label{Th1St2}
$\theta_\varepsilon<\pi$.
\end{step}

\proof[Proof of Step~\ref{Th1St2}]
Assume by contradiction that $\theta_\varepsilon=\pi$. Then it follows from \eqref{Th1Eq7} that 
$$\frac{d}{d\theta}\left.\left[Z\(x_\varepsilon,\exp_{x_\varepsilon}\(\theta v\),t_\varepsilon\)\right]\right|_{\theta=\pi}=\<\nabla_x\gamma\(y_\varepsilon,t_\varepsilon\),\nu_\varepsilon\(v\)\>-\partial_\theta\kappa\(\pi,t_\varepsilon\)=0$$
for all $v\in T_{x_\varepsilon}\S^n$ such that $\left|v\right|=1$, where $\nu_\varepsilon\(v\)=\frac{d}{d\theta}\left.\exp_{x_\varepsilon}\(\theta v\)\right|_{\theta=\pi}$. By observing that $\nu_\varepsilon\(-v\)=-\nu_\varepsilon\(v\)$, we then obtain a contradiction with $\partial_\theta\kappa\(\pi,t_\varepsilon\)>0$. This ends the proof of Step~\ref{Th1St2}.
\endproof

Remark that it follows from Steps~\ref{Th1St1} and~\ref{Th1St2} that for small $\varepsilon$, the function $Z$ is differentiable in a neighborhood of the point $\(x_\varepsilon,y_\varepsilon,t_\varepsilon\)$.

\smallskip
Our next result is as follows:

\begin{step}\label{Th1St3}
There exists a constant $\Lambda_0=\Lambda_0\(I,\phi,n\)>0$ such that 
\begin{multline}\label{Th1St3Eq1}
\partial_t\kappa\(\theta_\varepsilon,t_\varepsilon\)\le\frac{\Lambda_0^{-1}\partial^2_\theta\kappa\(\theta_\varepsilon,t_\varepsilon\)}{\(1+\partial_\theta\kappa\(\theta_\varepsilon,t_\varepsilon\)^2\)^{3/2}}\\
+\frac{\Lambda_0\(\theta_\varepsilon^{-1}\kappa\(\theta_\varepsilon,t_\varepsilon\)+\partial_\theta\kappa\(\theta_\varepsilon,t_\varepsilon\)\)\kappa\(\theta_\varepsilon,t_\varepsilon\)\partial_\theta\kappa\(\theta_\varepsilon,t_\varepsilon\)}{\(1+\partial_\theta\kappa\(\theta_\varepsilon,t_\varepsilon\)^2\)^{1/2}}
\end{multline}
for small $\varepsilon>0$.
\end{step}

\proof[Proof of Step~\ref{Th1St3}]
We let $\tau_{\varepsilon}:\[0,\theta_{\varepsilon}\]\to\S^n$ be a minimizing geodesic from $x_{\varepsilon}$ to $y_{\varepsilon}$. It follows from \eqref{Th1Eq7} that
$$\nabla_x Z\(x_\varepsilon,y_\varepsilon,t_\varepsilon\)=0,\quad\nabla_y Z\(x_\varepsilon,y_\varepsilon,t_\varepsilon\)=0,\quad\text{and}\quad\partial_tZ\(x_\varepsilon,y_\varepsilon,t_\varepsilon\)\ge0$$
which give
\begin{equation}\label{Th1St3Eq2}
\left\{\begin{aligned}
&\nabla_x \gamma\(x_\varepsilon,t_\varepsilon\)=\partial_\theta\kappa\(\theta_\varepsilon,t_\varepsilon\)\tau_\varepsilon'\(0\)\\
&\nabla_x\gamma\(y_\varepsilon,t_\varepsilon\)=\partial_\theta\kappa\(\theta_\varepsilon,t_\varepsilon\)\tau_\varepsilon'\(\theta_\varepsilon\)\\
&\partial_t\gamma\(y_\varepsilon,t_\varepsilon\)-\partial_t\gamma\(x_\varepsilon,t_\varepsilon\)\ge\partial_t\kappa\(\theta_\varepsilon,t_\varepsilon\).
\end{aligned}\right.
\end{equation}
By using \eqref{Th1Eq3} and \eqref{Th1St3Eq2}, we obtain
\begin{equation}\label{Th1St3Eq3}
\partial_t\kappa\(\theta_\varepsilon,t_\varepsilon\)\le\frac{A_\varepsilon}{\(1+\partial_\theta\kappa\(\theta_\varepsilon,t_\varepsilon\)^2\)^{3/2}}\\
+\frac{B_\varepsilon+\partial_\theta\kappa\(\theta_\varepsilon,t_\varepsilon\)^2C_\varepsilon}{\(1+\partial_\theta\kappa\(\theta_\varepsilon,t_\varepsilon\)^2\)^{1/2}},
\end{equation}
where
\begin{align*}
A_\varepsilon&:=\frac{\[\nabla_x^2\gamma\(y_\varepsilon,t_\varepsilon\)\]\(\tau_\varepsilon'\(\theta_\varepsilon\),\tau_\varepsilon'\(\theta_\varepsilon\)\)}{\psi\(\gamma\(y_\varepsilon,t_\varepsilon\)\)}-\frac{\[\nabla_x^2\gamma\(x_\varepsilon,t_\varepsilon\)\]\(\tau_\varepsilon'\(0\),\tau_\varepsilon'\(0\)\)}{\psi\(\gamma\(x_\varepsilon,t_\varepsilon\)\)},\allowdisplaybreaks\\
B_\varepsilon&:=\frac{\Delta_x\gamma\(y_\varepsilon,t_\varepsilon\)-\[\nabla_x^2\gamma\(y_\varepsilon,t_\varepsilon\)\]\(\tau_\varepsilon'\(\theta_\varepsilon\),\tau_\varepsilon'\(\theta_\varepsilon\)\)}{\psi\(\gamma\(y_\varepsilon,t_\varepsilon\)\)}\\
&\qquad-\frac{\Delta_x\gamma\(x_\varepsilon,t_\varepsilon\)-\[\nabla_x^2\gamma\(x_\varepsilon,t_\varepsilon\)\]\(\tau_\varepsilon'\(0\),\tau_\varepsilon'\(0\)\)}{\psi\(\gamma\(x_\varepsilon,t_\varepsilon\)\)},\allowdisplaybreaks\\
C_\varepsilon&:=n\(\frac{\psi'\(\gamma\(y_\varepsilon,t_\varepsilon\)\)}{\psi\(\gamma\(y_\varepsilon,t_\varepsilon\)\)^2}-\frac{\psi'\(\gamma\(x_\varepsilon,t_\varepsilon\)\)}{\psi\(\gamma\(x_\varepsilon,t_\varepsilon\)\)^2}\).
\end{align*}
Since $\kappa\(\theta_\varepsilon,t_\varepsilon\),\partial_\theta\kappa\(\theta_\varepsilon,t_\varepsilon\)>0$ and $\partial_\theta^2\kappa\(\theta_\varepsilon,t_\varepsilon\)<0$, in order to obtain \eqref{Th1St3Eq1}, it remains to prove that there exist constants $c_1,c_2,c_3>0$ depending only on $I$, $\phi$, and $n$ such that 
\begin{equation}\label{Th1St3Eq4}
\left\{\begin{aligned}
&A_\varepsilon\le c_1\partial_\theta^2\kappa\(\theta_\varepsilon,t_\varepsilon\)\\
&B_\varepsilon\le c_2\theta_\varepsilon^{-1}\kappa\(\theta_\varepsilon,t_\varepsilon\)^2\partial_\theta\kappa\(\theta_\varepsilon,t_\varepsilon\)\\
&C_\varepsilon\le c_3\kappa\(\theta_\varepsilon,t_\varepsilon\)
\end{aligned}\right.
\end{equation}
for small $\varepsilon$. We begin with proving the last estimate in \eqref{Th1St3Eq4}. Remark that by using Step~\ref{Th1St1}, we obtain  
$$\frac{\kappa\(\theta_\varepsilon,t_\varepsilon\)}{\varepsilon}=\frac{1}{\varepsilon}\int_0^{\theta_\varepsilon}\partial_\theta\kappa\(s,t_\varepsilon\)ds=\frac{1}{\varepsilon}\int_0^{\theta_\varepsilon}\frac{\overline\lambda e^{-\eta t_\varepsilon}}{\sqrt{\delta+s}}ds\ge\frac{\overline\lambda e^{-\eta t_\varepsilon}\theta_\varepsilon}{\varepsilon\sqrt{\delta+\theta_\varepsilon}}\to\infty$$
as $\varepsilon\to0$, which, together with \eqref{Th1Eq7}, implies that 
\begin{equation}\label{AddedEq}
\gamma\(y_\varepsilon,t_\varepsilon\)-\gamma\(x_\varepsilon,t_\varepsilon\)\le 2\kappa\(\theta_\varepsilon,t_\varepsilon\)
\end{equation}
for small $\varepsilon$. Since $\gamma\(D_T\)\subseteq\gamma_0\(\S^n\)$ and $\rho_0\(\S^n\)\subseteq I$, by applying the mean value theorem together with \eqref{AddedEq}, we obtain  
\begin{equation}\label{Th1St3Eq5}
C_\varepsilon\le2n\sup_{\gamma_0\(\S^n\)}\(\frac{\psi'}{\psi^2}\)'\kappa\(\theta_\varepsilon,t_\varepsilon\)\le2n\sup_I\(\frac{\phi'}{\phi}\)'\kappa\(\theta_\varepsilon,t_\varepsilon\)
\end{equation}
for small $\varepsilon$ which gives the last estimate in \eqref{Th1St3Eq4}. Now we prove the first two estimates in \eqref{Th1St3Eq4}. We let $\(v_{\varepsilon,1}\(0\),\dotsc,v_{\varepsilon,n}\(0\)\)$ be an orthonormal basis of $T_{x_\varepsilon}\S^n$ such that $v_{\varepsilon,n}\(0\)=\tau_\varepsilon'\(0\)$. For any $i\in\left\{1,\dotsc,n\right\}$, we let $\varphi_{\varepsilon,i}$ be a smooth function on $\[0,\theta_\varepsilon\]$ such that
\begin{equation}\label{Th1St3Eq6}
\varphi_{\varepsilon,i}\(0\)=\frac{1}{\sqrt{\psi\(\gamma\(x_\varepsilon,t_\varepsilon\)\)}}\quad\text{and}\quad\varphi_{\varepsilon,i}\(\theta_\varepsilon\)=\frac{\delta_i}{\sqrt{\psi\(\gamma\(y_\varepsilon,t_\varepsilon\)\)}}
\end{equation}
with $\delta_i:=1$ in case $i\ne n$ and $\delta_i:=-1$ in case $i=n$. For any $r\ge0$ and $\theta\in\[0,\theta_\varepsilon\]$, we define 
$$\tau_{\varepsilon,i}\(r,\theta\):=\exp_{\tau_\varepsilon\(\theta\)}\(r\varphi_{\varepsilon,i}\(\theta\)v_{\varepsilon,i}\(\theta\)\),$$
where $v_{\varepsilon,i}:\[0,\theta_\varepsilon\]\to T\S^n$ is the parallel transport of $v_{\varepsilon,i}\(0\)$ along $\tau_\varepsilon$. By using \eqref{Th1St3Eq6}, we obtain
\begin{equation}\label{Th1St3Eq7}
A_\varepsilon=\frac{d^2}{dr^2}\left.\[\gamma\(\tau_{\varepsilon,n}\(r,\theta_\varepsilon\),t_\varepsilon\)-\gamma\(\tau_{\varepsilon,n}\(r,0\),t_\varepsilon\)\]\right|_{r=0}
\end{equation}
and
\begin{equation}\label{Th1St3Eq8}
B_\varepsilon=\sum_{i=1}^{n-1}\frac{d^2}{dr^2}\left.\[\gamma\(\tau_{\varepsilon,i}\(r,\theta_\varepsilon\),t_\varepsilon\)-\gamma\(\tau_{\varepsilon,i}\(r,0\),t_\varepsilon\)\]\right|_{r=0}.
\end{equation}
On the other hand, for any $i\in\left\{1,\dotsc,n\right\}$, since $\partial_\theta\kappa\(\cdot,t_\varepsilon\)>0$ on $\(0,\pi\)$, it follows from \eqref{Th1Eq7} that
\begin{multline}\label{Th1St3Eq9}
\gamma\(\tau_{\varepsilon,i}\(r,\theta_\varepsilon\),t_\varepsilon\)-\gamma\(\tau_{\varepsilon,i}\(r,0\),t_\varepsilon\)-\varepsilon\le\kappa\(d\(\tau_{\varepsilon,i}\(r,0\),\tau_{\varepsilon,i}\(r,\theta_\varepsilon\)\),t_\varepsilon\)\\
\le\kappa\(\int_0^{\theta_\varepsilon}\left|\partial_\theta\tau_{\varepsilon,i}\(r,\theta\)\right|d\theta,t_\varepsilon\)
\end{multline}
for all $i\in\left\{1,\dotsc,n\right\}$ for small $r\ge0$ with equality in case $r=0$. Moreover, by using \eqref{Th1St3Eq2}, we obtain
\begin{align}\label{Th1St3Eq10}
&\frac{d}{dr}\left.\[\gamma\(\tau_{\varepsilon,i}\(r,\theta_\varepsilon\),t_\varepsilon\)-\gamma\(\tau_{\varepsilon,i}\(r,0\),t_\varepsilon\)\]\right|_{r=0}\nonumber\\
&\quad=\partial_\theta\kappa\(\theta_\varepsilon,t_\varepsilon\)\(\varphi_{\varepsilon,i}\(\theta_\varepsilon\)\<v_{\varepsilon,n}\(\theta_\varepsilon\),v_{\varepsilon,i}\(\theta_\varepsilon\)\>-\varphi_{\varepsilon,i}\(0\)\<v_{\varepsilon,n}\(0\),v_{\varepsilon,i}\(0\)\>\)\nonumber\allowdisplaybreaks\\
&\quad=\partial_\theta\kappa\(\theta_\varepsilon,t_\varepsilon\)\int_0^{\theta_\varepsilon}\<v_{\varepsilon,n}\(\theta\),\varphi_{\varepsilon,i}'\(\theta\)v_{\varepsilon,i}\(\theta\)\>d\theta\nonumber\allowdisplaybreaks\\
&\quad=\partial_\theta\kappa\(\theta_\varepsilon,t_\varepsilon\)\int_0^{\theta_\varepsilon}\<\partial_\theta\tau_{\varepsilon,i}\(0,\theta\),\partial_r\partial_\theta\tau_{\varepsilon,i}\(0,\theta\)\>d\theta\nonumber\\
&\quad=\frac{d}{dr}\left.\left[\kappa\(\int_0^{\theta_\varepsilon}\left|\partial_\theta\tau_{\varepsilon,i}\(r,\theta\)\right|d\theta,t_\varepsilon\)\right]\right|_{r=0}.
\end{align}
It follows from \eqref{Th1St3Eq9} and \eqref{Th1St3Eq10} that
\begin{align}
&\frac{d^2}{dr^2}\left.\[\gamma\(\tau_{\varepsilon,i}\(r,\theta_\varepsilon\),t_\varepsilon\)-\gamma\(\tau_{\varepsilon,i}\(r,0\),t_\varepsilon\)\]\right|_{r=0}\nonumber\\
&\quad\le\frac{d^2}{dr^2}\left.\left[\kappa\(\int_0^{\theta_\varepsilon}\left|\partial_\theta\tau_{\varepsilon,i}\(r,\theta\)\right|d\theta,t_\varepsilon\)\right]\right|_{r=0}\nonumber\allowdisplaybreaks\\
&\quad=\partial_\theta^2\kappa\(\theta_\varepsilon,t_\varepsilon\)\(\frac{d}{dr}\left.\left[\int_0^{\theta_\varepsilon}\left|\partial_\theta\tau_{\varepsilon,i}\(r,\theta\)\right|d\theta\right]\right|_{r=0}\)^2\nonumber\\
&\qquad+\partial_\theta\kappa\(\theta_\varepsilon,t_\varepsilon\)\frac{d^2}{dr^2}\left.\left[\int_0^{\theta_\varepsilon}\left|\partial_\theta\tau_{\varepsilon,i}\(r,\theta\)\right|d\theta\right]\right|_{r=0}.\label{Th1St3Eq11}
\end{align}
By proceeding as in \eqref{Th1St3Eq10} and using \eqref{Th1St3Eq6}, we obtain\addtolength{\textheight}{4pt}
\begin{align}\label{Th1St3Eq12}
&\frac{d}{dr}\left.\left[\int_0^{\theta_\varepsilon}\left|\partial_\theta\tau_{\varepsilon,i}\(r,\theta\)\right|d\theta\right]\right|_{r=0}\nonumber\\
&\quad=\varphi_{\varepsilon,i}\(\theta_\varepsilon\)\<v_{\varepsilon,n}\(\theta_\varepsilon\),v_{\varepsilon,i}\(\theta_\varepsilon\)\>-\varphi_{\varepsilon,i}\(0\)\<v_{\varepsilon,n}\(0\),v_{\varepsilon,i}\(0\)\>\nonumber\\
&\quad=\left\{\begin{aligned}&0&&\text{if }i\ne n\\&\frac{1}{\sqrt{\psi\(\gamma\(x_\varepsilon,t_\varepsilon\)\)}}+\frac{1}{\sqrt{\psi\(\gamma\(y_\varepsilon,t_\varepsilon\)\)}}&&\text{if }i=n.\end{aligned}\right.
\end{align}
Moreover, since $\gamma\(D_T\)\subseteq\gamma_0\(\S^n\)$ and $\rho_0\(\S^n\)\subseteq I$, we obtain
\begin{align}\label{Th1St3Eq13}
\frac{1}{\sqrt{\psi\(\gamma\(x_\varepsilon,t_\varepsilon\)\)}}+\frac{1}{\sqrt{\psi\(\gamma\(y_\varepsilon,t_\varepsilon\)\)}}\ge2\inf_{\gamma_0\(\S^n\)}\frac{1}{\sqrt{\psi}}\ge2\inf_I\frac{1}{\sqrt{\phi}}\,.
\end{align}
By differentiating twice, we obtain
\begin{align}\label{Th1St3Eq14}
&\frac{d^2}{dr^2}\left.\left[\int_0^{\theta_\varepsilon}\left|\partial_\theta\tau_{\varepsilon,i}\(r,\theta\)\right|d\theta\right]\right|_{r=0}=\int_0^{\theta_\varepsilon}\big(\left|\partial_\theta\tau_{\varepsilon,i}\(0,\theta\)\right|^{-1}\big(\left|\partial_r\partial_\theta\tau_{\varepsilon,i}\(0,\theta\)\right|^2\nonumber\\
&\qquad\quad+\<\partial_\theta\tau_{\varepsilon,i}\(0,\theta\),\partial_r^2\partial_\theta\tau_{\varepsilon,i}\(0,\theta\)\>\big)\nonumber\\
&\qquad\quad-\left|\partial_\theta\tau_{\varepsilon,i}\(0,\theta\)\right|^{-3}\<\partial_\theta\tau_{\varepsilon,i}\(0,\theta\),\partial_r\partial_\theta\tau_{\varepsilon,i}\(0,\theta\)\>^2\big)d\theta\nonumber\allowdisplaybreaks\\
&\qquad=\int_0^{\theta_\varepsilon}\big(\left|v_{\varepsilon,n}\(\theta\)\right|^{-1}\big(\varphi_{\varepsilon,i}'\(\theta\)^2\left|v_{\varepsilon,i}\(\theta\)\right|^2\nonumber\\
&\qquad\quad-\<v_{\varepsilon,n}\(\theta\),R\(\varphi_{\varepsilon,i}\(\theta\)v_{\varepsilon,i}\(\theta\),v_{\varepsilon,n}\(\theta\)\)\varphi_{\varepsilon,i}\(\theta\)v_{\varepsilon,i}\(\theta\)\>\big)\nonumber\\
&\qquad\quad-\left|v_{\varepsilon,n}\(\theta\)\right|^{-3}\<v_{\varepsilon,n}\(\theta\),\varphi_{\varepsilon,i}'\(\theta\)v_{\varepsilon,i}\(\theta\)\>^2\big)d\theta\nonumber\\
&\qquad=\left\{\begin{aligned}&\int_0^{\theta_\varepsilon}\(\varphi_{\varepsilon,i}'\(\theta\)^2-\varphi_{\varepsilon,i}\(\theta\)^2\)d\theta&&\text{if }i\ne n\\&0&&\text{if }i=n,\end{aligned}\right.
\end{align}
where $R$ is the curvature tensor of $\(\S^n,g_{\S^n}\)$. Since $\partial_\theta^2\kappa\(\theta_\varepsilon,t_\varepsilon\)\le0$, the first estimate in \eqref{Th1St3Eq4} follows from \eqref{Th1St3Eq7} and \eqref{Th1St3Eq11}--\eqref{Th1St3Eq14}. Now, we prove the second estimate in \eqref{Th1St3Eq4}. In case $i\ne n$, by integrating by parts, we obtain 
\begin{multline}\label{Th1St3Eq15}
\int_0^{\theta_\varepsilon}\(\varphi_{\varepsilon,i}'\(\theta\)^2-\varphi_{\varepsilon,i}\(\theta\)^2\)d\theta=\varphi_{\varepsilon,i}\(\theta_\varepsilon\)\varphi_{\varepsilon,i}'\(\theta_\varepsilon\)-\varphi_{\varepsilon,i}\(0\)\varphi_{\varepsilon,i}'\(0\)\\
-\int_0^{\theta_\varepsilon}\(\varphi_{\varepsilon,i}''\(\theta\)+\varphi_{\varepsilon,i}\(\theta\)\)\varphi_{\varepsilon,i}\(\theta\)d\theta.
\end{multline}
By using \eqref{Th1St3Eq15} with the function $\varphi_{\varepsilon,i}$ defined as
$$\varphi_{\varepsilon,i}\(\theta\):=\frac{1}{\sin\(\theta_\varepsilon\)}\(\frac{\sin\(\theta_\varepsilon-\theta\)}{\sqrt{\psi\(\gamma\(x_\varepsilon,t_\varepsilon\)\)}}+\frac{\sin\(\theta\)}{\sqrt{\psi\(\gamma\(y_\varepsilon,t_\varepsilon\)\)}}\)\quad\forall\theta\in\[0,\theta_\varepsilon\],$$
we obtain\addtolength{\textheight}{-4pt}
\begin{align}\label{Th1St3Eq16}
&\int_0^{\theta_\varepsilon}\(\varphi_{\varepsilon,i}'\(\theta\)^2-\varphi_{\varepsilon,i}\(\theta\)^2\)d\theta=\varphi_{\varepsilon,i}\(\theta_\varepsilon\)\varphi_{\varepsilon,i}'\(\theta_\varepsilon\)-\varphi_{\varepsilon,i}\(0\)\varphi_{\varepsilon,i}'\(0\)\nonumber\\
&\quad=\frac{1}{\sin\(\theta_\varepsilon\)}\bigg[\(\frac{1}{\psi\(\gamma\(x_\varepsilon,t_\varepsilon\)\)}+\frac{1}{\psi\(\gamma\(y_\varepsilon,t_\varepsilon\)\)}\)\cos\(\theta_\varepsilon\)\nonumber\\
&\qquad-\frac{2}{\sqrt{\psi\(\gamma\(x_\varepsilon,t_\varepsilon\)\)\psi\(\gamma\(y_\varepsilon,t_\varepsilon\)\)}}\bigg]\nonumber\allowdisplaybreaks\\
&\quad=\frac{\cos\(\theta_\varepsilon\)}{\sin\(\theta_\varepsilon\)}\bigg(\frac{1}{\sqrt{\psi\(\gamma\(x_\varepsilon,t_\varepsilon\)\)}}-\frac{1}{\sqrt{\psi\(\gamma\(y_\varepsilon,t_\varepsilon\)\)}}\bigg)^2\nonumber\\
&\qquad-\frac{2\tan\(\theta_\varepsilon/2\)}{\sqrt{\psi\(\gamma\(x_\varepsilon,t_\varepsilon\)\)\psi\(\gamma\(y_\varepsilon,t_\varepsilon\)\)}}\nonumber\\
&\quad\le\frac{1}{\theta_\varepsilon}\bigg(\frac{1}{\sqrt{\psi\(\gamma\(x_\varepsilon,t_\varepsilon\)\)}}-\frac{1}{\sqrt{\psi\(\gamma\(y_\varepsilon,t_\varepsilon\)\)}}\bigg)^2.
\end{align}
By proceeding as in \eqref{Th1St3Eq5}, we obtain 
\begin{equation}\label{Th1St3Eq17}
\left|\frac{1}{\sqrt{\psi\(\gamma\(x_\varepsilon,t_\varepsilon\)\)}}-\frac{1}{\sqrt{\psi\(\gamma\(y_\varepsilon,t_\varepsilon\)\)}}\right|\le\sup_I\(\frac{\phi'}{\sqrt\phi}\)\kappa\(\theta_\varepsilon,t_\varepsilon\)
\end{equation}
for small $\varepsilon$. By using \eqref{Th1St3Eq14}, \eqref{Th1St3Eq16}, and \eqref{Th1St3Eq17}, we obtain
\begin{equation}\label{Th1St3Eq18}
\sum_{i=1}^{n-1}\frac{d^2}{dr^2}\left.\left[\int_0^{\theta_\varepsilon}\left|\partial_\theta\tau_{\varepsilon,i}\(r,\theta\)\right|d\theta\right]\right|_{r=0}\le\sup_I\(\frac{\phi'}{\sqrt\phi}\)^2\frac{\kappa\(\theta_\varepsilon,t_\varepsilon\)^2}{\theta_\varepsilon}
\end{equation}
for small $\varepsilon$. The second estimate in \eqref{Th1St3Eq4} then follows from \eqref{Th1St3Eq8}, \eqref{Th1St3Eq11}, \eqref{Th1St3Eq12}, and \eqref{Th1St3Eq18}. This ends the proof of Step~\ref{Th1St3}.
\endproof

We can now end the proof of Theorem~\ref{Th1}.

\proof[End of proof of Theorem~\ref{Th1}] 
By applying Step~\ref{Th1St3} and observing that $\kappa\(\theta_\varepsilon,t_\varepsilon\)\le2\theta_\varepsilon\partial_\theta\kappa\(\theta_\varepsilon,t_\varepsilon\)$ and $\kappa\(\theta_\varepsilon,t_\varepsilon\)\partial_\theta\kappa\(\theta_\varepsilon,t_\varepsilon\)\le2\overline\lambda^2e^{-2\eta t_\varepsilon}$, we obtain
\begin{equation}\label{Th1Eq8}
-2\eta\big(\sqrt{\delta+\theta_\varepsilon}-\sqrt\delta\big)\le\frac{-\Lambda_0^{-1}}{2\big(\delta+\theta_\varepsilon+\overline\lambda^2e^{-2\eta t_\varepsilon}\big)^{3/2}}+\frac{6\Lambda_0\overline\lambda^2e^{-2\eta t_\varepsilon}}{\big(\delta+\theta_\varepsilon+\overline\lambda^2e^{-2\eta t_\varepsilon}\big)^{1/2}}
\end{equation}
Since $\delta<1$, $\theta_\varepsilon<\pi$, and $e^{-2\eta t_\varepsilon}\le1$, it follows from \eqref{Th1Eq8} that
$$1\le4\Lambda_0\big(1+\pi+\overline\lambda^2\big)\Big(\eta\sqrt{\(1+\pi\)\big(1+\pi+\overline\lambda^2\big)}+3\Lambda_0\overline\lambda^2\Big).$$
which gives a contradiction when $\overline\lambda$ and $\eta$ are smaller than some constants depending only on $I$, $\phi$, and $n$. This proves that for such values of $\overline\lambda$ and $\eta$, we have $Z\le0$ in $U_T$ and so
\begin{equation}\label{Th1Eq9}
\sup_{\S^n}\left|\nabla\gamma\(\cdot,t\)\right|\le\overline\lambda\delta^{-1/2}e^{-\eta t}\qquad\forall t\in\[0,T\].
\end{equation}
Since $\left|\nabla\rho\right|=\phi\(\rho\)\left|\nabla\gamma\right|$, it follows from \eqref{Th1Eq9} that 
\begin{equation}\label{Th1Eq10}
\sup_{\S^n}\left|\nabla\rho\(\cdot,t\)\right|\le\sup_I\(\phi\)\overline\lambda\delta^{-1/2}e^{-\eta t}\qquad\forall t\in\[0,T\].
\end{equation}
It then follows from classical theory of parabolic equations that $\rho\(\cdot,t\)$ exists for all $t\ge0$. Moreover, it follows from \eqref{Th1Eq10} that $\rho\(\cdot,t\)$ converges exponentially to a constant. This ends the proof of Theorem~\ref{Th1}.
\endproof

\section{Proof of Theorem~\ref{Th2}}\label{Sec3}

This section is devoted to the proof of Theorem~\ref{Th2}. We first prove the following result:

\begin{lemma}\label{Lem}
Assume that $n=1$, $0\in I$, $\phi$ is even, and $\phi''\(0\)>0$. Let $\psi$ and $\Gamma$ be as in \eqref{Th1Eq2} with $\overline\rho=0$, $J:=\Gamma\(I\)$, and for any $\tau>0$ and $k\in\N\backslash\left\{0\right\}$, $D_{1,k,\tau}:=\[0,\pi/\(2k\)\)\cup\(\pi/\(2k\),\pi/k\]\times\[0,\tau\)$ and $D_{2,k,\tau}:=\[0,\pi/k\]\times\[0,\tau\)$. Then for any $\sigma\in\(0,1/2\)$ and $\mu>0$, there exists $\varepsilon_0>0$ such that for any $\tau>0$ and $k\in\N$ such that $k^2\tau<\varepsilon_0$ and $k>1/\varepsilon_0$, there exist $\zeta_1\in C^\infty\(D_{1,k,\tau},J\)\cap C^0\(D_{2,k,\tau},J\)$ and $\zeta_2\in C^\infty\(D_{2,k,\tau},J\)$ such that
\begin{equation}\label{LemEq1}
\partial_t\zeta_i\le\frac{1}{\psi\(\zeta_i\)}\frac{\partial^2_\theta\zeta_i}{\(1+\(\partial_\theta\zeta_i\)^2\)^{3/2}}+\frac{\psi'\(\zeta_i\)}{\psi\(\zeta_i\)^2}\frac{\(\partial_\theta\zeta_i\)^2}{\sqrt{1+\(\partial_\theta\zeta_i\)^2}}\quad\text{in }D_{i,k,\tau}
\end{equation}
for $i\in\left\{1,2\right\}$ and the function $\zeta:=\max\(\zeta_1,\zeta_2\)$ is such that
\begin{enumerate}
\item[(A1)]$\displaystyle\zeta_1\(\pi/\(2k\),t\)<\zeta_2\(\pi/\(2k\),t\)(=\zeta\(\pi/\(2k\),t\))\quad\forall t\in\[0,\tau\)$,\vspace{5pt}
\item[(A2)]$\displaystyle\left|\zeta\(\theta,0\)-\zeta\(\theta',0\)\right|<\mu\left|\theta-\theta'\right|^\sigma\quad\forall\theta,\theta'\in\[0,\pi/k\],\vspace{5pt}$
\item[(A3)]$\displaystyle\zeta\(0,t\)=\zeta\(\pi/k,t\)=0\quad\forall t\in\[0,\tau\),$\vspace{5pt}
\item[(A4)]$\partial_\theta\zeta\(0,t\)\to\infty$ and $\partial_\theta\zeta\(\pi/k,t\)\to-\infty$ as $t\to\tau$.
\end{enumerate}
\end{lemma}

\proof[Proof of Lemma~\ref{Lem}]
We fix $p\in\(2/\(1-\sigma\),4\)$. We let $\tau>0$ and $k\in\N\backslash\left\{0\right\}$ to be chosen later on so that $k$ is large and $k^2\tau$ is small. For any $\(\theta,t\)\in D_{2,k,\tau}$, we define
$$\zeta_1\(\theta,t\):=\min\(\frac{c_1\theta}{\(\(\tau-t\)^p+\theta^2\)^{1/p}}\,,\frac{c_1\(\pi/k-\theta\)}{\[\(\tau-t\)^p+\(\pi/k-\theta\)^2\]^{1/p}}\)$$
and
$$\zeta_2\(\theta,t\):=c_1A_k\(\sin\(k\theta\)-c_2k^2t\),$$
where $A_k:=2^{2/p}\(\pi/k\)^{1-2/p}$ and $c_1$ and $c_2$ are positive constants independent of $\theta$, $t$, $k$, and $\tau$ to be fixed later on. Note that $1-2/p>\sigma$. It is easy to check that $\zeta_1\in C^\infty\(D_{1,k,\tau},J\)\cap C^0\(D_{2,k,\tau},J\)$, $\zeta_2\in C^\infty\(D_{2,k,\tau},J\)$, and (A2)--(A4) hold true for small $\tau$ and large $k$. If moreover $k^2\tau$ is small, then we obtain that (A1) holds true. It remains to prove that \eqref{LemEq1} holds true. Since $\phi$ is even and $\phi''\(0\)>0$, we obtain that $\psi$ is also even and $\psi''\(0\)>0$. By applying the mean value theorem and since $\psi'\(0\)=0$ and $\zeta_1\(\theta,t\)\ge0$, we obtain
\begin{equation}\label{LemEq2}
\frac{\psi'\(\zeta_1\(\theta,t\)\)}{\psi\(\zeta_1\(\theta,t\)\)^2}\ge\inf_{\zeta_1\(D_{2,k,\tau}\)}\(\frac{\psi'}{\psi^2}\)'\zeta_1\(\theta,t\)
\end{equation}
for all $\(\theta,t\)\in D_{2,k,\tau}$. Moreover, direct calculations give
\begin{align}
&\partial_t\zeta_1\(\theta,t\)=\frac{c_1\theta\(\tau-t\)^{p-1}}{\(\(\tau-t\)^p+\theta^2\)^{1+1/p}}\le\frac{c_1\theta}{\(\(\tau-t\)^p+\theta^2\)^{2/p}}\,,\label{LemEq3}\allowdisplaybreaks\\
&\frac{\partial^2_\theta\zeta_1\(\theta,t\)}{\(1+\(\partial_\theta\zeta_1\(\theta,t\)\)^2\)^{3/2}}\nonumber\\
&\qquad=-\frac{2c_1\theta\(3\(\tau-t\)^p+\(1-2/p\)\theta^2\)\(\(\tau-t\)^p+\theta^2\)^{1+2/p}}{p\big(\(\(\tau-t\)^p+\theta^2\)^{2+2/p}+c_1^2\(\(\tau-t\)^p+\(1-2/p\)\theta^2\)^2\big)^{3/2}}\nonumber\\
&\qquad\ge-\frac{6\theta}{c_1^2p\(1-2/p\)^3\(\(\tau-t\)^p+\theta^2\)^{1-2/p}}\,,\label{LemEq4}
\end{align}
and
\begin{align}\label{LemEq5}
&\frac{\zeta_1\(\theta,t\)\(\partial_\theta\zeta_1\(\theta,t\)\)^2}{\sqrt{1+\(\partial_\theta\zeta_1\(\theta,t\)\)^2}}\nonumber\\
&\qquad=\frac{c_1^3\theta\(\(\tau-t\)^p+\(1-2/p\)\theta^2\)^2\(\(\tau-t\)^p+\theta^2\)^{-1-2/p}}{\big(\(\(\tau-t\)^p+\theta^2\)^{2+2/p}+c_1^2\(\(\tau-t\)^p+\(1-2/p\)\theta^2\)^2\big)^{1/2}}\nonumber\\
&\qquad\ge\frac{c_1^3\(1-2/p\)^2\theta}{\sqrt{\(\tau^p+\(\pi/\(2k\)\)^2\)^{2/p}+c_1^2}\(\(\tau-t\)^p+\theta^2\)^{2/p}}
\end{align}
for all $\(\theta,t\)\in\[0,\pi/\(2k\)\)\times\[0,\tau\)$. Since $p<4$, it follows from \eqref{LemEq2}--\eqref{LemEq5} that
\begin{multline}\label{LemEq9}
\frac{1}{\psi\(\zeta_1\(\theta,t\)\)}\frac{\partial^2_\theta\zeta_1\(\theta,t\)}{\(1+\(\partial_\theta\zeta_1\(\theta,t\)\)^2\)^{3/2}}+\frac{\psi'\(\zeta_1\(\theta,t\)\)}{\psi\(\zeta_1\(\theta,t\)\)^2}\frac{\(\partial_\theta\zeta_1\(\theta,t\)\)^2}{\sqrt{1+\(\partial_\theta\zeta_1\(\theta,t\)\)^2}}\\
-\partial_t\zeta_1\(\theta,t\)\ge\frac{\theta}{\(\(\tau-t\)^p+\theta^2\)^{2/p}}\bigg(-\frac{6\(\tau^2+\(\pi/\(2k\)\)^2\)^{-1+4/p}}{c_1^2p\(1-2/p\)^3}\hspace{-6pt}\sup_{\zeta_1(D_{2,k,\tau})}\frac{1}{\psi}\\
+\frac{c_1^3\(1-2/p\)^2}{\sqrt{\(\tau^p+\(\pi/\(2k\)\)^2\)^{2/p}+c_1^2}}\inf_{\zeta_1(D_{2,k,\tau})}\(\frac{\psi'}{\psi^2}\)'-c_1\bigg)
\end{multline}
for all $\(\theta,t\)\in\[0,\pi/\(2k\)\)\times\[0,\tau\)$ provided
$$\inf_{\zeta_1(D_{2,k,\tau})}\(\frac{\psi'}{\psi^2}\)'\ge0.$$
Moreover, direct calculations give
\begin{equation}\label{LemEq10}
\zeta_1\(D_{2,k,\tau}\)=\big[0,c_1\(\pi/k\)^{1-2/p}\big).
\end{equation}
Since $2<p<4$, we obtain
\begin{equation}\label{LemEq10a}
1-\frac{2}{p}>0\qquad\text{and}\qquad-1+\frac{4}{p}>0.
\end{equation}
Since $\psi'\(0\)=0$, and $\psi''\(0\)>0$, it follows from \eqref{LemEq9}--\eqref{LemEq10a} that \eqref{LemEq1} holds true for $i=1$ for small $\tau$ and large $k$ provided the constant $c_1$ is chosen large enough so that $c_1>\(1-2/p\)^{-2}\psi\(0\)^2/\psi''\(0\)$. With regard to the function $\zeta_2$, we obtain
\begin{align}
\partial_t\zeta_2\(\theta,t\)&=-c_1c_2k^2A_k,\label{LemEq11}\\
\frac{\partial^2_\theta\zeta_2\(\theta,t\)}{\(1+\(\partial_\theta\zeta_2\(\theta,t\)\)^2\)^{3/2}}&=-\frac{c_1k^2A_k\sin\(k\theta\)}{\(1+c_1^2k^2A_k^2\cos\(k\theta\)^2\)^{3/2}}\nonumber\\
&\ge-c_1k^2A_k\,,\label{LemEq12}
\end{align}
and
\begin{equation}\label{LemEq13}
\frac{\(\partial_\theta\zeta_2\(\theta,t\)\)^2}{\sqrt{1+\(\partial_\theta\zeta_2\(\theta,t\)\)^2}}=\frac{c_1^2k^2A_k^2\cos\(k\theta\)^2}{\sqrt{1+c_1^2k^2A_k^2\cos\(k\theta\)^2}}\in\[0,c_1^2k^2A_k^2\]
\end{equation}
for all $\(\theta,t\)\in D_{2,k,\tau}$. It follows from \eqref{LemEq11}--\eqref{LemEq13} that
\begin{multline}\label{LemEq14}
\frac{1}{\psi\(\zeta_2\(\theta,t\)\)}\frac{\partial^2_\theta\zeta_2\(\theta,t\)}{\(1+\(\partial_\theta\zeta_2\(\theta,t\)\)^2\)^{3/2}}+\frac{\psi'\(\zeta_2\(\theta,t\)\)}{\psi\(\zeta_2\(\theta,t\)\)^2}\frac{\(\partial_\theta\zeta_2\(\theta,t\)\)^2}{\sqrt{1+\(\partial_\theta\zeta_2\(\theta,t\)\)^2}}\\
-\partial_t\zeta_2\(\theta,t\)\ge c_1k^2A_k\bigg(-\sup_{\zeta_2(D_{2,k,\tau})}\frac{1}{\psi}+A_kc_1\min\(\inf_{\zeta_2(D_{2,k,\tau})}\frac{\psi'}{\psi^2},0\)+c_2\bigg)
\end{multline}
for all $\(\theta,t\)\in D_{2,k,\tau}$. Moreover, direct calculations give
\begin{equation}\label{LemEq15}
\zeta_2\(D_{2,k,\tau}\)=\big(-c_1c_2A_kk^2\tau ,c_1A_k\big].
\end{equation}
It follows from \eqref{LemEq10a} and \eqref{LemEq15} that for every $\varepsilon>0$, if $k^2\tau<\varepsilon$ and $k>1/\varepsilon$, then
\begin{equation}\label{LemEq16}
\zeta_2\(D_{2,k,\tau}\)\subset\(-2^{2/p}\pi^{1-2/p}\varepsilon^{2-2/p} c_1c_2,2^{2/p}\pi^{1-2/p}\varepsilon^{1-2/p} c_1\).
\end{equation}
By continuity of $1/\psi$ and $\psi'/\psi^2$ and since $\psi\(0\)>0$ , $A_k\to0$ as $k\to\infty$ and $\varepsilon^{1-2/p}\to0$ as $\varepsilon\to0$, it follows from \eqref{LemEq14} and \eqref{LemEq16} that if the constant $c_2$ is chosen so that $c_2>1/\psi\(0\)$, then there exists $\varepsilon_0>0$ such that 
\begin{equation}\label{LemEq17}
-\sup_{\zeta_2(D_{2,k,\tau})}\frac{1}{\psi}+A_kc_1\min\(\inf_{\zeta_2(D_{2,k,\tau})}\frac{\psi'}{\psi^2},0\)+c_2>0
\end{equation}
for all $\tau>0$ and $k\in\N$ such that $k^2\tau<\varepsilon_0$ and $k>1/\varepsilon_0$. By putting together \eqref{LemEq14} and \eqref{LemEq17}, we obtain that \eqref{LemEq1} holds true with $i=2$. This ends the proof of Lemma~\ref{Lem}. 
\endproof

Now we can prove Theorem~\ref{Th2}.

\proof[Proof of Theorem~\ref{Th2}]
We fix $\sigma<1/2$, $\lambda>0$, and we define
\begin{equation}\label{Th2Eq2}
\mu:=\lambda\inf_I\frac{1}{\phi}\,.
\end{equation}
We let $\psi$ and $\Gamma$ be as in \eqref{Th1Eq2} and $J$, $\tau$, $k$, $D_{1,k,\tau}$, $D_{2,k,\tau}$, $\zeta_1$, $\zeta_2$, and $\zeta$ be as in Lemma~\ref{Lem}. By using (A1)--(A3) and since $\zeta_1\in C^\infty\(D_{1,k,\tau},J\)\cap C^0\(D_{2,k,\tau},J\)$ and $\zeta_2\in C^\infty\(D_{2,k,\tau},J\)$, we obtain that there exists $\varepsilon_0,a_0,b_0\in\R$ such that 
\begin{equation}\label{Th2Eq2b}
b_0<\min\(\frac{a_0\pi}{2k},\(\mu-\varepsilon_0\)\(\frac{\pi}{2k}\)^\sigma,\sup J\)
\end{equation}
and
\begin{equation}\label{Th2Eq2a}
\zeta\(\theta,0\)<\widetilde\gamma_0\(\theta\)\qquad\forall\theta\in\(0,\pi/k\),
\end{equation}
where
$$\widetilde\gamma_0\(\theta\):=\left\{\begin{aligned}&\min\(a_0\theta,\(\mu-\varepsilon_0\)\theta^\sigma,b_0\)&&\text{if }0\le\theta\le\frac{\pi}{2k}\\&\widetilde\gamma_0\(\frac{\pi}{k}-\theta\)&&\text{if }\frac{\pi}{2k}<\theta\le\frac{\pi}{k}\,.\end{aligned}\right.$$
For any $\varepsilon>0$ and $\theta\in\[0,\pi/k\]$, we then define
$$\widetilde\gamma_0^{\(\varepsilon\)}\(\theta\):=\left\{\begin{aligned}&f_\varepsilon\(a_0\theta,f_\varepsilon\(\(\mu-\varepsilon_0\)\theta^\sigma,b_0\)\)&&\text{if }0\le\theta\le\frac{\pi}{2k}\\&\widetilde\gamma_0^{\(\varepsilon\)}\(\frac{\pi}{k}-\theta\)&&\text{if }\frac{\pi}{2k}<\theta\le\frac{\pi}{k}\,,\end{aligned}\right.$$
where
$$f_\varepsilon\(\xi_1,\xi_2\):=\frac{1}{2}\[\xi_1+\xi_2-\varepsilon\,\eta\(\frac{\xi_2-\xi_1}{\varepsilon}\)\]\qquad\forall\xi_1,\xi_2\in\R$$
and $\eta:\R\to\(0,\infty\)$ is a smooth, even cutoff function such that $\eta\(\theta\)=\theta$ for all $\theta\in\[1,\infty\)$ and $\eta'\(\theta\)>0$ for all $\theta\in\(0,1\)$. By using \eqref{Th2Eq2a}, it is easy to see that for small $\varepsilon$, $\widetilde\gamma_0^{\(\varepsilon\)}\in C^\infty\(\[0,\pi/k\],J\)$ and $\widetilde\gamma_0^{\(\varepsilon\)}\to\widetilde\gamma_0$ in $C^{0,1}\(\[0,\pi/k\]\)$. Hence, by using \eqref{Th2Eq2b} and remarking that $\widetilde\gamma_0^{\(\varepsilon\)}\(\theta\)=\widetilde\gamma_0\(\theta\)=a_0\theta$ for small $\theta$, we obtain that for small $\varepsilon$, $\widetilde\gamma_0^{\(\varepsilon\)}$ is such that
\begin{enumerate}
\item[(B1)]$\displaystyle\zeta\(\theta,0\)\le\widetilde\gamma_0^{\(\varepsilon\)}\(\theta\)\quad\forall\theta\in\[0,\pi/k\],$\vspace{5pt}
\item[(B2)]$\displaystyle\big|\widetilde\gamma_0^{\(\varepsilon\)}\(\theta\)-\widetilde\gamma_0^{\(\varepsilon\)}\(\theta'\)\big|<\mu\left|\theta-\theta'\right|^\sigma\quad\forall\theta,\theta'\in\[0,\pi/k\],$\vspace{5pt}
\item[(B3)]$\widetilde\gamma_0^{\(\varepsilon\)}\(0\)=\widetilde\gamma_0^{\(\varepsilon\)}\(\pi/k\)=\widetilde\gamma_0^{\(\varepsilon\)}{''}\(0\)=\widetilde\gamma_0^{\(\varepsilon\)}{''}\(\pi/k\)=0.$
\end{enumerate}
In what follows, we fix $\varepsilon$ small enough so that (B1)--(B3) hold true. Since $\widetilde\gamma_0^{\(\varepsilon\)}\in C^\infty\(\[0,\pi/k\],J\)$, the classical theory of parabolic equations (see for instance Lieberman~\cite{Lie}*{Theorem~8.2}) gives the existence of a solution $\widetilde\gamma\in C^\infty\(\[0,\pi/k\]\times\[0,T\)\)$ of the problem
\begin{equation}\label{Th2Eq3}
\left\{\begin{aligned}&\hspace{-3pt}\partial_t\widetilde\gamma=\frac{1}{\psi\(\widetilde\gamma\)}\frac{\partial^2_\theta\widetilde\gamma}{\(1+\(\partial_\theta\widetilde\gamma\)^2\)^{3/2}}+\frac{\psi'\(\widetilde\gamma\)}{\psi\(\widetilde\gamma\)^2}\frac{\(\partial_\theta\widetilde\gamma\)^2}{\sqrt{1+\(\partial_\theta\widetilde\gamma\)^2}}\hspace{-6pt}&&\text{in }[0,\pi/k]\times[0,T)\\&\hspace{-3pt}\widetilde\gamma\(\cdot,0\)=\widetilde\gamma_0^{\(\varepsilon\)}&&\text{on }[0,\pi/k]\\&\hspace{-3pt}\widetilde\gamma\(0,\cdot\)=\widetilde\gamma\(\pi/k,\cdot\)=0&&\text{on }[0,T),\end{aligned}\right.
\end{equation}
where $T\in\(0,\infty\]$ is the maximal existence time for $\widetilde\gamma$. Moreover, since $\widetilde\gamma_0^{\(\varepsilon\)}\(\[0,\pi/k\]\)\subseteq J$, it follows from the maximum principle that $\widetilde\gamma\(\[0,\pi/k\]\times\[0,T\)\)\subseteq J$. By using (A1)  and \eqref{LemEq1} and integrating by parts, we obtain that $\zeta$ is a weak subsolution of the equation in \eqref{Th2Eq3}, i.e.
$$\int_0^{\tau'}\hspace{-5pt}\int_0^{\pi/k}\bigg(\eta\partial_t\zeta+\frac{1}{\psi\(\zeta\)}\frac{\partial_\theta\zeta\partial_\theta\eta}{\sqrt{1+\(\partial_\theta\zeta\)^2}}-2\frac{\psi'\(\zeta\)}{\psi\(\zeta\)^2}\frac{\(\partial_\theta\zeta\)^2\eta}{\sqrt{1+\(\partial_\theta\zeta\)^2}}\bigg)d\theta dt\le0$$
for all $\tau'\in\(0,\tau\)$ and $\eta\in C^1\(D_{2,k,\tau}\)$ such that $\eta\ge0$ in $D_{2,k,\tau}$ and $\eta\(0,\cdot\)=\eta\(\pi/k,\cdot\)=0$ on $[0,\tau)$. We define $\omega:=\zeta-\widetilde\gamma$. It follows from (A3), (B1), and \eqref{Th2Eq3} that $\omega\le0$ on $\left\{0,\pi/k\right\}\times\[0,\min\(T,\tau\)\)$ and $\[0,\pi/k\]\times\left\{0\right\}$. By applying the mean value theorem, we obtain that for any $\tau'\in\(0,\min\(T,\tau\)\)$, there exist $a_1,a_2,b_1,b_2\in L^\infty\(D_{2,k,\tau'}\)$ such that $\inf\left\{a_1\(\theta,t\):\,\(\theta,t\)\in D_{2,k,\tau'}\right\}>0$ and
\begin{equation}\label{Th2Eq4}
\int_0^{\tau'}\hspace{-5pt}\int_0^{\pi/k}\hspace{-5pt}\(\eta\partial_t\omega+\(a_1\partial_\theta\omega+a_2\omega\)\partial_\theta\eta+\(b_1\partial_\theta\omega+b_2\omega\)\eta\)d\theta dt\le0
\end{equation}
for all $\eta\in C^1\(D_{2,k,\tau'}\)$ such that $\eta\ge0$ in $D_{2,k,\tau'}$ and $\eta\(0,\cdot\)=\eta\(\pi/k,\cdot\)=0$ on $[0,\tau')$. By applying a weak comparison principle (see for instance Lieberman~\cite{Lie}*{Corollary~6.16}), it follows from \eqref{Th2Eq4} that
\begin{equation}\label{Th2Eq5}
\zeta\(\theta,t\)\le\widetilde\gamma\(\theta,t\)
\end{equation}
for all $\(\theta,t\)\in D_{2,k,\min\(T,\tau\)}$. It follows from (A3), (A4), and \eqref{Th2Eq5} that $T\le\tau$. Note that by using similar arguments as in \eqref{Th2Eq4}--\eqref{Th2Eq5}, we obtain that $\widetilde\gamma$ is the unique solution of \eqref{Th2Eq3}. It then follows from classical theory of parabolic equations that 
\begin{equation}\label{Th2Eq6}
\lim_{t\to T}\sup_{D_{2,k,t}}\left|\partial_\theta\widetilde\gamma\right|=\infty.
\end{equation}
Indeed, if \eqref{Th2Eq6} is not true, then $T=\infty$ (see for instance Lieberman~\cite{Lie}*{Theorems~8.3 and~12.1}), which is in contradiction with $T\le\tau$. We let $\gamma:\S^1\times\[0,T\)\to\R$ be the function defined as
$$\gamma\(\(\cos\theta,\sin\theta\),t\):=\left\{\begin{aligned}&\widetilde\gamma\(\theta-j\pi/k,t\)&&\text{if }j\text{ is even}\\&-\widetilde\gamma\(\(j+1\)\pi/k-\theta,t\)&&\text{if }j\text{ is odd}\end{aligned}\right.$$
for all $\(\theta,t\)\in\[j\pi/k,\(j+1\)\pi/k\)\times\[0,\min\(T,\tau\)\)$, $j\in\left\{0,\dotsc,2k-1\right\}$. Since $\psi'\(0\)=0$, it follows from \eqref{Th2Eq3} and (B3) that $\partial^2_\theta\widetilde\gamma\(j\pi/k,t\)=0$ for all $t\in\[0,T\)$ and $j\in\left\{0,\dotsc,2k-1\right\}$ which implies that $\gamma$ is a smooth solution of \eqref{Th1Eq3}. By using (B2), \eqref{Th2Eq2}, \eqref{Th2Eq6}, and the change of functions \eqref{Th1Eq2}, we then obtain the existence of $\rho_0\in C^\infty\(\S^n,I\)$ such that \eqref{Eq4} holds true, the solution of \eqref{Eq1} exists and $\partial_x\rho\(\cdot,t\)$ blows up as $t\to T$. This ends the proof of Theorem~\ref{Th2}.
\endproof

\end{document}